\documentclass[12pt]{article}
\usepackage[cp1251]{inputenc}
\begin{document}
\large

\begin{center}
{\Large \bf Solvability of linear non-homogeneous higher order
differential equation in the Schwartz space}
\end{center}

\vskip5mm

\begin{center}
\noindent Valerii Hr. Samoilenko, Yuliia I. Samoilenko\\
Department of Mathematical Physics, \\
Taras Shevchenko National
University of Kyiv, Kyiv, Ukraine, 01601 \\
valsamyul@gmail.com, vsam@univ.kiev.ua, yusam@univ.kiev.ua
\end{center}

\vskip5mm

\begin{minipage}{150mm}
{\bf Abstract.} There is studied problem on solvability of linear
non-homogeneous differential equation of higher even order. There is
proved the theorem on necessary and sufficient conditions on
existence of solutions to the equation in the Schwartz space.
\end{minipage}

\vskip10mm


\vskip10mm


\vskip10mm

\section{Introduction and statement of problem} One of the most
important problem in the theory of differential equations of
mathematical physics is the one on the existence of solutions to
various differential equations. In particular, while studying
soliton solutions there is arisen a problem on existence of solution
in the Schwartz space. For example, the problem on existence of
these solutions to the following equation with the Schr{\"{o}}dinger
operator
\begin{equation} \label{eq_Shrodinger}
- \frac{d^2 v}{dx^2} + q v = f, \quad x \in{\mathbf R},
\end{equation}
has been arisen while constructing asymptotic soliton like solutions
to singularly perturbed Korteweg-de Vries equation \cite{Sam1}.

The problem was solved through theory of pseudodifferential
operators in \cite{SamShrod} where there were found necessary and
sufficient conditions for existence of solutions to equation
(\ref{eq_Shrodinger}) in the Schwartz space.

On the other hand, while studying soliton solutions to the higher
Korteweg-de Vries equations and KdV-like equations there is appeared
the same problem for similar differential equations with the
operator of higher order. So, there is come up problem on studying
necessary and sufficient conditions on existence of the solutions to
the following equation
\begin{equation} \label{high_eq}
L v = f, \quad x \in{\mathbf R},
\end{equation}
in the Schwartz space, where differential operator $L$ is written as
\begin{equation} \label{operator_L}
L = - \sum\limits_{m=1}^n a_m \frac{d^{2m}}{d x^{2m}} + q
\end{equation}
with constant coefficients and $ q $ being a function.

\section{Main result} Let $ \mathrm{S} ({\mathbf R}) $ be the Schwartz space.

The main result of the paper is proposition.

{\bf Theorem.} {\it Let the following conditions be fulfilled:

$1^0.$ the coefficients $ a_1, a_2, \dots , a_n $ are nonnegative
constants and $ a_n > 0 $;

$2^0.$ $ q(x) = q_0 + q_1(x) $, where constant $q_0 <0$ and the
function $ q_1(x) \in \mathrm{S}({\mathbf R}) $. 


If the kernel of the operator $ L : \mathrm{S}({\mathbf R}) \to
\mathrm{S}({\mathbf R}) $ is trivial, then equation (\ref{high_eq})
has a solution in the space $ \mathrm{S}({\mathbf R}) $ for any
function $ f \in \mathrm{S}({\mathbf R}) $.

Otherwise, if the kernel of the operator $ L: \mathrm{S}({\mathbf
R}) \to \mathrm{S}({\mathbf R}) $ is not trivial, then equation
(\ref{high_eq}) has a solution in the space $ \mathrm{S}({\mathbf
R}) $ if and only if the function $ f \in \mathrm{S}({\mathbf R}) $
satisfies the condition of orthogonality in the form
\begin{equation}\label{ort_cond}
\int\limits_{-\infty}^{+\infty} f (x) v_0(x) d x = 0
\end{equation}
for any $ v_0 \in ker~L $.}

\section{Necessary definitions and statements} Let us remind some
notions and results from theory of pseudodifferential operators.
These ones are using below under proving the theorem.

For any function $ h \in \mathrm{S} ({\mathbf R}) $ there is defined
the Fourier transform as
$$
F[h](\xi) = \int\limits_{-\infty}^{+\infty} e^{- i \xi x} h(x) d x.
$$

Due to properties of the Fourier transform it's possible to define
the inverse Fourier transform as
\begin{equation} \label{symbol}
p \left( x, \frac{d}{dx} \right) h (x) = \frac{1}{2\pi}
\int\limits_{-\infty}^{+\infty} e^{i x \xi} \, p(x, \xi) F[h](\xi)
\,   d \xi
\end{equation}
for any differential operator
$$
p \left( x, \frac{d}{dx} \right) = \sum\limits_{k = 0}^n a_k(x)
\frac{d^k}{dx^k} , \quad x \in {\mathbf R}.
$$

Here
$$ p(x, \xi) = \sum\limits_{k = 0}^n a_k(x) (- i\xi)^k , \quad x,
\xi \in {\mathbf R},
$$
is called a symbol of the differential operator $ p \left( x,
\frac{d}{dx} \right) $.

Let $ S^m $ be a set of such symbols $ p(x, \xi) $ $ \in $ $
\mathrm{C}^{\infty} ({\mathbf R}^2) $ that for any $ k $, $ l \in
{\mathbf N}\cup\{0\} $ the inequality
$$
\left| p_{(l)}^{(k)} (x, \xi) \right| \le C_{k l} \left( 1 +
|\xi|\right) ^{m - k}, \quad (x, \xi) \in {\mathbf R}^2,
$$
is satisfied, where
$$
p_{(l)}^{(k)} (x, \xi) = \frac{\partial^{k + l}}{\partial\xi^k
\partial x^l} \, p(x, \xi), \quad (x, \xi) \in {\mathbf R}^2,
$$
and $ C_{kl} $ are some constants \cite{Hor}.

By $ S^m_0 $ denote a set of such symbols $ p(x, \xi) \in S^m $,
that
$$
|p(x, \xi) | \le M(x) \left(1 + |\xi| \right)^m,
$$
where value $ M(x) \to 0 $ as $ |x| \to + \infty $.

Let $ \mathrm{H}_s({\mathbf R}) $, $ s \in {\mathbf R} $, be a
Sobolev space \cite{GR1}, i.e. a space of such generalized functions
$ g \in \mathrm{S}^* ({\mathbf R}) $ that their Fourier transform $
F[g](\xi) $ satisfies condition
\begin{equation}\label{norma}
|| g ||_s^2 =  \int\limits_{-\infty}^{+\infty} (1 + |\xi|^2)^s \, |
F[g](\xi) |^2 \, d \xi  < \infty .
\end{equation}

\section{Proof of the main result} Proving the theorem contains two
steps. Firstly, we show that the operator $ L: \mathrm{H}_{s+2n}
({\mathbf R}) \rightarrow \mathrm{H}_s ({\mathbf R})$ of form
(\ref{operator_L}) is the Noether operator for any $ s \in {\mathbf
R} $. Later we prove that the solution to equation (\ref{high_eq})
belongs to the space $ \mathrm{S}({\mathbf R}) $.

Let us consider symbol of the differential operator $ L $ having a
form
\begin{equation} \label{symbol_L}
p (x, \xi) = -\sum\limits_{m=1}^{n} a_m \xi^{2m} + q(x).
\end{equation}

It's obviously that $ p(x, \xi) $ belongs to the set $ S^{2n} $
because of inequality
$$
\left| \frac{\partial^{k+l}}{\partial\xi^k \partial x^l} p(x, \xi)
\right| \le C_{k l} (1 +|\xi|)^{2n-k}, \quad k, l \in {\mathbf N}
\cup\{0\}.
$$
Moreover,
$$
\frac{\partial^{\,l}}{\partial x^l} \, p(x, \xi) \in S_0^{2n}, \quad
l \in{\mathbf N}.
$$

Accordingly to assumptions of the theorem the operator 
$ L: \mathrm{H}_{s+2n} ({\mathbf R}) \rightarrow \mathrm{H}_{s}
({\mathbf R}) $ satisfies all conditions of theorem 3.4 \cite{GR1}
for any $ s \in {\mathbf R} $. So, it is the Noether operator. As
consequence the operator $ L: \mathrm{H}_{s + 2n} ({\mathbf R}) \to
\mathrm{H}_{s} ({\mathbf R}) $ is normally solvable.

Thus, if kernel of the operator $ L^* $ is nontrivial then for
solvability of differential equation (\ref{high_eq}) in the space $
\mathrm{H}_{s} ({\mathbf R}) $ it's necessary and sufficient that
condition of orthogonality
\begin{equation} \label{adjoint}
<f, ker L^*> = 0
\end{equation}
is satisfied.

Since
$$
L^* = - \sum\limits_{m=1}^n a_m \frac{d^{2m}}{d x^{2m}} + q(x)
$$
then $ ker (L^*) \subset \bigcap\limits_{s \in {\mathbf R}}
\mathrm{H}_s({\mathbf R})$ \cite{GR1}.

Using Sobolev embedding theorems for the spaces $ \mathrm{H}_{s}
({\mathbf R}) $, $ s \in {\mathbf R} $, we have $ v_0^* \in \bar
{\mathrm{C}}_0^{\infty}({\mathbf R}) $ for any element $ v_0^* \in
ker {L^*} $.

Under the orthogonality condition (\ref{adjoint}) taking into
account theorem 3.4 \cite{GR1} we deduce that the solution $ v(x) $
of equation (\ref{high_eq}) is a such that $ v \in \bigcap\limits_{s
\in {\mathbf R}} \mathrm{H}_s({\mathbf R})$.

Arguing as above we obtain $ v \in \bar {\mathrm{C}}_0^{\infty}
({\mathbf R}). $

Now let us show that moreover $ v \in \mathrm{S}({\mathbf R}) $.
Indeed, since the function $ v \in \bar {\mathrm{C}}_0^{\infty}
({\mathbf R}) $ and the one satisfies equation
\begin{equation}\label{quikly_decres_fun}
- \sum\limits_{m=1}^n a_m \frac{d^{2m}}{d x^{2m}} v = - q v + f,
\end{equation}
where the function $ - q v + f \in \mathrm{S}({\mathbf R}) $, then
due to properties of elliptic pseudodifferential operators with
polynomial coefficients \cite{GR2}, we obtain that any solution to
equation (\ref{quikly_decres_fun}) from the space $
\mathrm{S}^*({\mathbf R}) $ belongs to the space $
\mathrm{S}({\mathbf R}) $. Thus, $ v \in \mathrm{S}({\mathbf R}) $.

Continuing this line of reasoning we see that $ v_0^* \in
\mathrm{S}({\mathbf R}) $. Kind of the operator $ L^* $ implies $
v_0^* = v_0 $. It means that orthogonality condition (\ref{adjoint})
is equivalent to the one (\ref{ort_cond}).

From the above consideration it follows that if the kernel of the
operator $ L : \mathrm{S}({\mathbf R}) \to \mathrm{S}({\mathbf R}) $
is trivial, i.e. the homogeneous equation $ L v = 0 $ has the only a
zero solution in the space $ \mathrm{S}({\mathbf R}) $, then
equation (\ref{high_eq}) has a solution in the space $
\mathrm{S}({\mathbf R}) $ under conditions a) and b) of theorem 3.4
\cite{GR1}.

The theorem is proved.

\section*{Conclusions}
Theorem on existence of a solution to the linear inhomogeneous
differential equation in the Schwartz space is proved. The theorem
can be used while studying asymptotic soliton like solutions to
singularly perturbed higher order partial differential equations of
integrable type with variable coefficients.


\begin{thebibliography}{99}

\bibitem{Sam1} {V.Hr.~Samoylenko and Yu.I.~Samoylenko, “Asymptotic
expansions for one-phase soliton-like solutions to the Korteweg–de
Vries equation with variable coefficients,” Ukranian Mathematical
Journal, 57 (2005), no.~1, 111~--~124.}

\bibitem{SamShrod} {V.Hr.~Samoylenko and Yu.I.~Samoylenko, "Existence of a
solution to the inhomogeneous equation with the one-dimensional
Schrodinger operator in the space of quickly decreasing functions"
Journal of Mathematical Sciences, 187 (2012), no.~1, 70~--~76.}

\bibitem{Hor} {L.~H\"{o}rmander, The Analysis of Linear Partial Differential
Operators: Pseudodifferential Operators, Springer, Berlin, 1985.}

\bibitem{GR1} {V.V.~Grushin, “Pseudodifferential operators in $R^n$ with bounded
symbols,” Functional Analysis and Its Applications July, 4 (1970),
no.~3, 202~--~212.}

\bibitem{GR2} {V.V.~Grushin, “On a class of elliptic pseudodifferential
operators degenerating on a submanifold,” Math. USSR-Sb., 13 (1971),
no.~2, 155~--~185.}

\end{thebibliography}
\end{document}